\documentclass[12pt]{amsart}
\usepackage{amsmath}
\usepackage{amstext}
\usepackage{amsfonts}
\usepackage{verbatim}
\usepackage{amssymb}

\renewcommand{\a}{\alpha}

\newcommand{\G}{\Gamma}
\renewcommand{\l}{\lambda}
\renewcommand{\i}{\infty}
\newcommand{\hh}{H^2(\G^2)}
\newcommand{\hho}{H^2(\G)}

\def\p{\varphi}

\newenvironment{keeptogether}{\pagebreak[0]\samepage}{}

\newtheorem{Theorem}{Theorem}[section]
\newtheorem{Proposition}[Theorem]{Proposition}
\newtheorem{Lemma}[Theorem]{Lemma}
\newtheorem{Corollary}[Theorem]{Corollary}
\newtheorem{Conjecture}[Theorem]{Conjecture}
\newtheorem{Remark}[Theorem]{Remark}
\newtheorem{Question}[Theorem]{Question}
\newtheorem{Definition}[Theorem]{Definition}

\theoremstyle{definition}
\newtheorem{Example}[Theorem]{Example}

\newenvironment{theorem}{\begin{keeptogether}\begin{Theorem}\slshape}%
     {\end{Theorem}\end{keeptogether}}
\newenvironment{proposition}{\begin{keeptogether}\begin{Proposition}\slshape
}%
     {\end{Proposition}\end{keeptogether}}
\newenvironment{lemma}{\begin{keeptogether}\begin{Lemma}\slshape}%
{\end{Lemma}\end{keeptogether}}
\newenvironment{corollary}{\begin{keeptogether}\begin{Corollary}\slshape}%
     {\end{Corollary}\end{keeptogether}}
     {\end{Conjecture}\end{keeptogether}}
     {\end{Remark}\end{keeptogether}}
     {\end{Question}\end{keeptogether}}
     {\end{Example}\end{keeptogether}}

     {\end{Definition}\end{keeptogether}}

\newcounter{subtheoremc}



\renewcommand{\em}{\sl}

\begin{document}

\title{$N_{\p}$-type quotient modules on the torus}

\author[K. Izuchi]{Keiji Izuchi}
\address{Department of Mathematics, Niigata University, 
Niigata, 950-2181, Japan}
\email{izuchi@math.sc.niigata-u.ac.jp}
\thanks{The first author is partially supported by Grant-in-Aid for
Scientific Research (No.16340037),
Ministry of Education, Science and Culture.}

\author[R. Yang]{Rongwei Yang }
\address{Department of Mathematics and Statistics, SUNY at Albany, 
Albany, NY 12047, U.S.A.}
\email{ryang@math.albany.edu}
\subjclass{Primary 46E20; Secondary 47A13}

\begin{abstract}
Structure of the quotient modules in $\hh$ is very complicated. A good understanding of some special
examples will shed light on the general picture. This paper studies the so-call $N_{\p}$-type quotient modules, 
namely, quotient modules of the form $\hh\ominus [z-\p]$, where $\p (w)$ is a function in the classical Hardy space $H^2(\G)$
and $[z-\p]$ is the submodule generated by $z-\p (w)$. This type of quotient modules serve as good examples in many studies. 
A notable feature of the $N_{\p}$-type quotient module is its close connections with some classical single variable 
operator theories.

\end{abstract}

\maketitle
\section{Introduction}
Let $H^2(\G^2)$ be the Hardy space on the two dimensional torus $\G^2$.
We denote by $z$ and $w$ the coordinate functions. Shift operators $T_z$ and $T_w$ on $H^2(\G^2)$ are defined by 
$T_zf=zf$ and $T_wf=wf$ for $f\in H^2(\G^2)$. Clearly, both $T_z$ and $T_w$ have infinite multiplicity.
A closed subspace $M$ of $H^2(\G^2)$ is called a {\em submodule} (over the algebra $H^{\infty}(D^2)$), if 
it is invariant under multiplications by functions in $H^{\infty}(D^2)$, where $D$ stands for the unit disk. Equivalently, M is a submodule 
if it is invariant for both $T_z$ and $T_w$. The quotient space $N:=H^2(\G^2)\ominus M$ is
called a {\em quotient module}. Clearly $T^*_zN\subset N$ and $T^*_wN\subset N$. And for this reason
$N$ is also said to be backward shift invariant. In the study here, it is necessary to distinguish the classical Hardy 
space in the variable $z$ and that in the variable $w$, for which we denote by $H^2(\G_z)$ and $H^2(\G_w)$, respectively.
$H^2(\G_z)$ and $H^2(\G_w)$ are thus different subspaces in $\hh$. We will simply write $H^2(\G)$ when there is no need 
to tell the difference. In $H^2(\G)$, it is well known as the Beurling theorem 
that if $M\subset H^2(\G)$ is invariant for $T_z$, then $M=qH^2(\G)$ for an inner function $q(z)$.
The structure of submodules in $H^2(\G^2)$ is much more
complex, and there is a great amount of works on this subject in recent years. A good reference of this work 
can be found in \cite{CG}. One natural approach to the problem is to find and study some relatively simple submodules,
and hope that the study will generate concepts and general techniques that will lead to
a better understanding of the general picture. This in fact has become an interesting and encouraging work.

In this paper, we look at submodules of the form $[z-\p(w)]$,
where  $\p$ is a function in $H^2(\G_w)$ with $\p\not=0$ and $[z-\p(w)]$ is the closure of $(z-\p)H^{\infty}(\G^2)$ in $H^2(\G^2)$. 
For simplicity we denote $[z-\p(w)]$ by $M_{\p}$.
One good way of studying $M_{\p}$ is through the so-called {\em two variable Jordan block} $(S_z,S_w)$ defined on the quotient module
\[N_{\p}:=H^2(\G^2)\ominus M_{\p}.\] For every 
quotient module $N$, the two variable Jordan block $(S_{z},S_{w})$ is  
the compression of the pair $(T_{z},T_{w})$ to $N$, or
more precisely, 
\[S_{z}f=P_Nzf,\ \ S_{w}f=P_Nwf,\ \ \ f\in N,\]
where $P_N:\hh \to N$ is the orthogonal projection. This paper studies interconnections between the quotient module
$N_{\p}$, the two variable Jordan block $(S_{z},S_{w})$ and the function $\p$.
Some related work has been done in \cite{IY, Y2, Y4}. By \cite{IY}, $N_\p\not=\{0\}$ if and only if $\p(D)\cap D\not=\emptyset$. 
For convenience, we let \[\Omega_{\p}=\{w\in D: |\p(w)|<1\},\] and assume throughout the paper that
$N_\p\not=\{0\}$, i.e., $\p(D)\cap D\not=\emptyset$. The paper is organized as follows.\\

Section 1 is introduction.\\

Section 2 introduces some useful tools and states a few related known results.\\

Section 3 studies the spectral properties of the operators $S_z$ and $S_w$. It is interesting to see how these 
properties depend on the function $\p$.\\

A notable phenomenon in many cases is the compactness of the defect operators $I-S_zS_z^*$ and $I-S^*_zS_z$. 
Section 4 aims to study how the compactness is related to the properties of $\p$.\\

The quotient module $N_{\p}$ has very rich structure. Indeed, when $\p$ is inner, $N_{\p}$ can be identified with
the tensor product of two well-known classical spaces, namely the quotient space $H^2(\G)\ominus \p H^2(\G)$ and the 
Bergman space $L^2_a(D)$. Section 5 makes a detailed study of this case.\\  

{\bf Acknowledgement.} This paper was finished when the second author was visiting the Niigata University. The hospitality and 
conveniences provided by its Department of Mathematics is greatly appreciated.

\section{Preliminaries}

For every $\lambda\in D$, we define a {\em left evaluation} 
operator $L(\lambda)$ from $\hh$ to $H^2(\G_w)$ and a 
{\em right evaluation} operator $R(\l)$ from $\hh$ to $H^2(\G_z)$ by
\[L(\lambda) f(w)=f(\lambda, w),\ \ \ R(\l)f(z)=f(z,\ \lambda),\ \ \ f\in \hh.\]
Clearly, $L(\l)$ and $R(\l)$ are operator-valued analytic functions over $D$.
Restrictions of $L(\l)$ and $R(\l)$ to quotient spaces $N$, $M\ominus zM$ and $M\ominus wM$ 
play key roles in the study here. The following lemma is from \cite{DY}.
\begin{lemma}
The restriction of $R(\l)$ to $M\ominus wM$ is equivalent to the characteristic operator function
for $S_w$.
\end{lemma}

The following spectral relations are thus clear. Details can be found in \cite{DY} and \cite{SF}.

(a) $\lambda \in \sigma (S_w)$ if and only if $R(\lambda): {M\ominus wM}\rightarrow  H^2(\G_z)$ is not invertible,

(b) $dim\,ker\,(S_w-\lambda I)=dim\, ker\,(R(\lambda)|_{M\ominus wM})$. 

(c) $S_w-\lambda I$ has a closed range if and only if $R(\lambda)(M\ominus wM)$ is closed,

(d) $S_w-\lambda I$ is Fredholm if and only if $R(\lambda)|_{M\ominus wM}$ is Fredholm, and in this case
\[ind\,(S_w-\lambda I)=ind\,(R(\lambda)|_{M\ominus wM}).\]

Restrictions $T^*_z|_{M\ominus zM}$ and $T^*_w|_{M\ominus wM}$ are also important here, and for simplicity they
are denoted by $D_{z}$ and $D_{w}$, respectively. Clearly,
\[D_{z}f (z,w)=\frac{f(z,w)-f(0,w)}{z},\ \ D_{w}f (z,w)=\frac{f(z,w)-f(z,0)}{w}.\]
And it is not hard to check that the ranges of $D_z$ and $D_w$ are subspaces of N. The following lemma (cf. \cite{Y2})
gives a description of the defect operators for $S_z$, and it will be used often.

\begin{lemma}
On a quotient module $N$,
\begin{itemize}
\item[(i)]
$S_z^*S_z+D_zD_z^*=I$;
\item[(ii)]
$S_zS_z^*+(L(0)|_N)^*L(0)|_N=I$.
\end{itemize}
\end{lemma}
A parallel version of Lemma 2.2 for $S_w$ will also be used. \\

The operator $D_z$ is a useful tool in this study. We first note that 
\[D^*_zf=P_M zf,\ \ \ f\in N.\]
So if $D^*_zf=0$, then $zf\in N$. Clearly $zf\in ker \,L(0)|_N$. 
Conversely, if $h$ is in $ker \,L(0)|_N$, then we can write 
$h=zh_0$. One checks easily that $h_0\in ker \,D_z^*.$ 
This observation shows that 
\[z\,ker\,D_z^*=ker \,L(0)|_N.    \]
So on $N_\p$, since $L(0)$ is injective (cf. \cite{IY}), $D_z^*$ has trivial kernel, i.e., the range $R(D_z)$ is dense in $N_\p$. 
The following theorem describes $R(D_z)$ in detail.

\begin{theorem}
Let $N$ be a quotient module of $H^2(\Gamma^2)$ and 
$M = H^2(\Gamma^2) \ominus N$.
Suppose that $R(D_z)$ is dense in $N$. Let $f \in N$.
Then $f \in R(D_z)$ if and only if there 
exists a positive constant 
$C_f$ depending on $f$ such that $|\langle S^*_zh, f\rangle | 
\le C_f \|L(0) h\|$ for every $h \in N$.
\end{theorem}

\begin{proof}
Suppose that $f \in R(D_z)$.
Let $g \in M \ominus zM$ with $T^*_z g = f$.
We have $g = zf + L(0)g$.
Then for $h \in N$,
\begin{eqnarray*}
|\langle S^*_zh, f\rangle | &=& |\langle h, zf\rangle| \\
&=& |\langle h, g - L(0)g\rangle|\\
&=& |\langle h, L(0)g\rangle|\\
&=& | \langle L(0) h, L(0)g\rangle|\\
&\le& \|L(0)g\| \|L(0) h\|.
\end{eqnarray*}

To prove the converse, suppose that there exists a positive 
constant $C_f$ satisfying 
\[|\langle S^*_zh, f\rangle | \le C_f \|L(0) h\|\] for every $h \in N$.
Since $L(0)$ on $N$ is one to one, we have a map $\Lambda$ defined by
$$
\Lambda: L(0) N \ni u(w) \to L(0)^{-1}u \to \langle S^*_z 
L(0)^{-1}u, f\rangle\in{\mathbb C}.
$$
Note that $L^{-1}_0u\in N$. 
Obviously, $\Lambda$ is linear and 
$$
|\Lambda u| = |\langle S^*_z L(0)^{-1}u, f\rangle| \le C_f \|L(0)
L(0)^{-1}u\| = C_f \|u\|.
$$
Hence by the Hahn-Banach theorem, $\Lambda$ is extendable to a 
bounded linear functional on $H^2(\Gamma_w)$ and there exists 
$v(w) \in H^2(\Gamma_w)$ satisfying
$\langle u, v \rangle = \Lambda u$ for every $u \in L(0) N$.
We have
$$
\langle u, v \rangle = \langle S^*_zL(0)^{-1}u, f\rangle
= \langle L(0)^{-1}u, zf\rangle.
$$
Since $v(w) \in H^2(\Gamma_w)$, $\langle u, v\rangle = 
\langle L(0)^{-1}u, v\rangle$.
Therefore
$$
\langle L(0)^{-1}u, zf-v\rangle = 0
$$
for every $u \in L(0) N$.
Since $L^{-1}_0(L(0)N) = N$, we get $zf-v \perp N$.
Hence $zf-v \in M$.
Since $v(w) \in H^2(\Gamma_w)$,
we have $T^*_z(zf-v) = f \in N$.
This implies that $zf-v \in M \ominus zM$.
Thus we get $f \in R(D_z)$.
\end{proof}

In the case of $N_\p$, \cite{IY} provides a very useful description of the functions in the space.
Let $\varphi(w)\in H^2(\Gamma_w)$.
For $f(w)\in H^2(\G_w)$, we formally define a function
$$
(T^*_\p f)(w)=\sum^\i_{n=0}a_nw^n,
$$
where
$$
a_n=\int^{2\pi}_0 \overline{\p}(e^{i\theta})f(e^{i\theta})e^{-in\theta}
d\theta/2\pi=\langle f(w), \p(w)w^n\rangle.
$$
Generally, $T^*_\p f$ may not be in $H^2(\G_w)$. When $T^*_\p f\in H^2(\G_w)$, 
we can define $T^{*2}_\p f=T^*_\p(T^*_\p f)$.
Inductively if $T^{*n}_\p f\in H^2(\G_w)$, 
we can define $T^{*(n+1)}_\p f=T^*_\p(T^{*n}_\p f)$.
For convenience, we let 
\[A_{\p}f(z,w)= \sum^\infty_{n=0}z^nT^{*n}_{\p}f(w)\]
be an operator defined at every $f\in H^2(\G_w)$ for which $A_{\p}f\in H^2(\Gamma^2)$.
Then it is shown in \cite{IY} that $L(0)$ is one-to-one on $N_\varphi$ and
$$
N_\varphi=\Big\{A_{\p}f: f\in H^2(\Gamma_w),\sum^\infty_{n=0}\|T^{*n}_{\p}f\|^2<\infty\Big\}.\leqno{(2.1)}
$$
It is easy to see that $L(0)A_\p f=f$. 
Moreover by \cite[Corollary 2.8]{IY},
$L(0)N_\p$ is dense in $H^2(\G_w)$.

The following two lemmas are needed for the study of $\sigma(S_z)$.

\begin{lemma}\label{lemma5.3}
Let $\p(w), g(w)\in H^2(\G_w)$ 
and $\psi(w)\in H^\i(\G_w)$.
Then $T^*_\p T^*_\psi g = T^*_{\psi\p}g.$
Moreover if $T^*_\p g\in H^2(\G_w)$, then 
$T^*_\psi T^*_\p g = T^*_{\psi\p}g$.
\end{lemma}

\begin{proof}
Let $n\ge0$. Then by the definitions above,
\[\langle T^*_\varphi T^*_\psi g, z^n\rangle = 
\langle g, \varphi\psi z^n\rangle=\langle T^*_{\varphi\psi} g, z^n\rangle.\]
Thus $T^*_\p T^*_\psi g=T^*_{\p\psi}g$.
Suppose that $T^*_\varphi g\in H^2(\Gamma_w)$.
We have $\overline{\varphi}g -T^*_\varphi g\in \overline{zH^1}$.
Hence 
\begin{eqnarray*}
\langle T^*_\psi T^*_\varphi g, z^n\rangle&=&
\langle T^*_\varphi g, \psi z^n\rangle\\
&=&\int^{2\pi}_0\overline{\p}(e^{i\theta})g(e^{i\theta})
\overline{\psi}(e^{i\theta})e^{-in\theta}d\theta/2\pi\\
&=&\langle g, \psi\varphi z^n\rangle.
\end{eqnarray*}
Thus we get our assertion.
\end{proof}

Let $w_0\in \Omega_\p$. The following lemma follows easily from the calculation
$$
T^*_\p\frac{1}{1-\overline{w}_0w}=\frac{\overline{\p(w_0)}}{1-\overline{w}_0w}.
$$

\begin{lemma}
For $w_0\in \Omega_\p$,we have
$$
\frac{1}{(1-\overline{\p(w_0)}z)(1-\overline{w}_0w)}\in N_\p.
$$
\end{lemma}

\section{the spectra of $S_z$ and $S_w$}

The spectra of $S_z$ and $S_w$ on $N_\p$ is evidently dependent on $\p$. 
This section aims to figure out how they 
are exactly related. Lemma 2.1 and the description in (2.1) are 
helpful to this end.

\begin{proposition}
$\overline{\varphi(D)\cap D}\subset \sigma(S_z)\subset 
\overline{\varphi(D)}\cap \overline{D}$.
\end{proposition}

\begin{proof}Let $w_0\in \varphi(D)\cap D$.
Then
\begin{eqnarray*}
S^*_z\Big(\frac{1}{(1-\overline{\p(w_0)}z)(1-\overline{w}_0w)}\Big)
&=& \sum^\infty_{n=1}\Big(\overline{\varphi(w_0)}^n(1-\overline{w}_0w)^{-1}
\Big)z^{n-1}\\
&=&\overline{\varphi(w_0)}\Big(
\frac{1}{(1-\overline{\p(w_0)}z)(1-\overline{w}_0w)}\Big).
\end{eqnarray*}
By Lemma 2.5, $\overline{\varphi(w_0)}$ is a point spectrum of $S^*_z$.
Thus we get $\overline{\varphi(D)\cap D} \subset\sigma(S_z)$.

Let $\lambda\notin \overline{\varphi(D)}$. 
Then $1/(\varphi(w)-\lambda)\in H^\infty(\Gamma_w)$.
Let $F\in N_\varphi$. We have
\begin{eqnarray*}
S^*_{1/(\varphi-\lambda)}F&=&S^*_{1/(\varphi-\lambda)}
\sum^\infty_{n=0}(T^{*n}_\varphi L(0)F)z^{n}\\
&=&\sum^\infty_{n=0}(T^{*n}_\varphi T^*_{1/(\varphi-\lambda)}
L(0)F)z^{n}\qquad\text{by Lemma 2.4.}
\end{eqnarray*}
Hence
\begin{eqnarray*}
S^*_{1/(\varphi-\lambda)}S^*_{z-\lambda}F&=&
\sum^\infty_{n=0}(T^{*n}_\varphi T^*_{1/(\varphi-\lambda)}L(0)
S^*_{z-\lambda}F)z^{n}\\
&=&\sum^\infty_{n=0}(T^{*n}_\varphi T^*_{1/(\varphi-\lambda)}
T^*_{\varphi-\lambda}L(0)F)z^{n}\\
&=&\sum^\infty_{n=0}(T^{*n}_\varphi L(0)F)z^{n}
\qquad\text{by Lemma 2.4}\\
&=&F.
\end{eqnarray*}
Also we have
\begin{eqnarray*}
S^*_{z-\lambda}S^*_{1/(\varphi-\lambda)}F&=&
\sum^\infty_{n=1}(T^{*n}_\varphi T^*_{1/(\varphi-\lambda)}L(0)F)z^{n-1}
-\overline{\lambda}\sum^\infty_{n=0}(T^{*n}_\varphi 
T^*_{1/(\varphi-\lambda)}L(0)F)z^{n}\\
&=&\sum^\infty_{n=1}(T^{*n}_\varphi T^*_\varphi T^*_{1/(\varphi-\lambda)}L(0)F)z^{n-1}-\overline{\lambda}\sum^\infty_{n=0}(T^{*n}_\varphi 
T^*_{1/(\varphi-\lambda)}L(0)F)z^{n}\\
&=&\sum^\infty_{n=0}(T^{*n}_\varphi T^*_{(\varphi-\lambda)}
T^*_{1/(\varphi-\lambda)}L(0)F)z^{n}\\
&=&F.
\end{eqnarray*}
Thus $(S_z-\lambda)^{-1}=S_{1/(\varphi-\lambda)}$
and hence $\lambda\notin\sigma(S_z)$.

Since $\|S_z\|\le1$, we have our assertion.
\end{proof}

For a submodule $M$ in $H^2(\G^2)$, the quotient space 
$M\ominus zM$ is a wandering subspace for the 
multiplication by $z$ and we have
\begin{equation*}
M=\sum^{\infty}_{n=0}\oplus\,z^n(M\ominus zM).
\end{equation*}
For a fixed $\l\in D$ and every $f\in M$, we write $f=\sum_{j=0}^{\infty}z^jf_j$ for some unique sequence
$\{f_j\}$ in $M\ominus zM$. So
\[f=\sum_{j=0}^{\infty}\l^jf_j+\sum_{j=0}^{\infty}(z^j-\l^j)f_j,\]
which means that $f=h_1+(z-\l)h_2$ for some $h_1\in M\ominus zM$ and $h_2\in M$.
If $h_1+(z-\l)h_2=0$, then $h_1+zh_2=\l h_2$, and hence
$|\l|^2\|h_2\|^2=\|h_1\|^2+\|h_2\|^2$,
which is possible only if $h_1=h_2=0$.
This observation shows that $M$ can be expressed as the direct sum
\begin{equation*}
 M=(M\ominus zM)+(z-\l)M,\tag{3.1}
\end{equation*}

We now look at the spectral properties of $S_w$.
\begin{proposition}
On $N_{\p}$,
\begin{itemize}
\item[(i)]
$\overline{\Omega_{\p}}\subset \sigma(S_w)$.
\item[(ii)]
$S_w-\a I$ is Fredholm for every $\a\in \Omega_{\p}$ and $ind\,(S_w-\a I)=-1$.
\end{itemize}
\end{proposition}
\begin{proof}

We use Lemma 2.1 to this end.

(i) It is sufficient to show $\Omega_{\p}\subset \sigma(S_w)$. If $\a\in \Omega_{\p}$, then for any function 
$(z-\p)h(z,w)$ in $M_\p\ominus wM_\p$, 
$(z-\p (\a))h(z,\a)$ vanishes at $\p(\a)$, and therefore 
$R(\a)({M_\p\ominus wM_\p})\subset (z-\p (\a))H^2(\G_z)\neq H^2(\G_z)$.
By Lemma 2.1, $\a\in \sigma(S_w)$.

(ii) It is equivalent to show that $R(\a)|_{M_\p\ominus wM_\p}$ is Fredholm with index $-1$.
We first show that $R(\a)$ is injective on ${M_\p\ominus wM_\p}$ for every $\a\in \Omega_{\p}$. Let $(z-\p)h(z,w)$ be in $M_\p$. Then
there is a sequence of polynomials $\{p_n(z,w)\}$ such that $(z-\p)p_n$ converges to $(z-\p)h$ in the norm of $\hh$.  
Since $R(\a)$ is a bounded operator, $(z-\p(\a))p_n(z,\a)$ converges to $(z-\p(\a))h(z,\a)$, which, by the fact
$|\p(\a)|<1$, implies that $p_n(z,\a)$ converges 
to $h(z,\a)$ in $H^2(\G_z)$.
Since for every $f \in H^2(\G_z)$, we have $\|\p f\|=\|\p\|\|f\|$ and hence 
\begin{equation*}
\|(z-\p)f\|\leq \|zf\|+\|\p f\|=(1+\|\p\|)\|f\|<\infty, \tag{3.2}
\end{equation*}
$(z-\p)p_n(z,\a)$ converges to  $(z-\p)h(z,\a)$ in $M_\p$. It follows that 
\[\lim_{n\to \infty}  (z-\p)\frac{p_n-p_n(\cdot,\a)}{w-\a} =(z-\p)\frac{h-h(\cdot,\a)}{w-\a},\]
which concludes that $(z-\p)\frac{h-h(\cdot,\a)}{w-\a}\in M_\p$.
If $(z-\p)h(z,w)$ is in $M_\p\ominus wM_\p$ such that $(z-\p(\a))h(z,\a)=0$, then $h(z,\a)=0$, and it follows from the observation above
that \[(z-\p)h=(w-\a)(z-\p)\frac{h}{w-\a}\in (w-\a)M_\p,\] and hence by (3.1) $(z-\p)h(z,w)=0$ which concludes that 
$R(\a)$ is injective on ${M_\p\ominus wM_\p}$. 

In the proof of (i), we showed that $R(\a)({M_\p\ominus wM_\p})\subset (z-\p (\a))H^2(\G_z)$. On the other hand, 
for every $g\in H^2(\G_z)$, $(z-\p)g$ is in $M_\p$ by (3.2), and by (3.1)  
\[(z-\p (\a))g\in R(\a)(M_\p)=R(\a)({M_\p\ominus wM_\p}).\] 
This shows that \[R(\a)({M_\p\ominus wM_\p})=(z-\p (\a))H^2(\G_z),\]
i.e., $R(\a)|_{M_\p\ominus wM_\p}$ has a closed range with codimension $1$, and this completes the proof in view of Lemma 2.1.
\end{proof}

\begin{corollary}
If $\p$ is bounded with $\|\p\|_{\infty}\leq 1$, then $\sigma(S_w)=\bar{D}$ and $\sigma_e(S_w)=\G$.
\end{corollary}
\begin{proof}
By Proposition 3.2 and the fact that $S_w$ is a contraction, $\sigma(S_w)=\bar{D}$ and $\sigma_e(S_w)\subset \G$. 
Since $ind(S_w)=-1$, $\sigma_e(S_w)$ is a closed curve, and therefore $\sigma_e(S_w)=\G$.
\end{proof}

We will mention another somewhat deeper consequence of Proposition 3.2 near the end of this section. Here we continue to study
the Fredholmness of $S_z$. Unfortunately, the techniques 
used for Proposition 3.2(ii) can not be applied directly to the case
here and a technical difficulty seems hard to overcome. So instead we use (3.1) in the case here. We begin with some simple 
observations.

\begin{lemma}\label{lemma5.1}
Let $\p(w)=b(w)h(w)$ be the inner-outer factorization of $\p(w)$.
Then $ker\,S^*_z = H^2(\G_w)\ominus b(w)H^2(\G_w)$.
\end{lemma}
\begin{proof}
Since the functions in $H^2(\G_w)\ominus b(w)H^2(\G_w)$ depend only on $w$, the inclusion 
\[H^2(\G_w)\ominus b(w)H^2(\G_w)\subset ker\,S^*_z \]
is easy to check.

If $f$ is a function in $N_{\p}$ such that $S^*_zf=0$, then $\bar{z}f$ is orthogonal to $H^2(\G^2)$ which means
$f$ is independent of the variable $z$. Since for every non-negative integer $j$
\[0=\langle(z-\p)w^j, f\rangle=\langle-\p w^j, f\rangle,\]
$f$ is in $H^2(\G_w)\ominus b(w)H^2(\G_w)$.
\end{proof}

\begin{theorem}\label{theorem5.4}
Let $\p(w)=b(w)h(w)$ be the inner-outer factorization of $\p$ and
\[\a=\displaystyle{\inf_{w\in D}}|h(w)|.\]
Then $S^*_z$ has a closed range if and only if $\a\not=0$, and
in this case $S^*_zN_\p=N_\p$.
\end{theorem}

\begin{proof}
Write $K_b= H^2(\Gamma_w)\ominus b(w)H^2(\Gamma_w)$.
By Lemma \ref{lemma5.1}, $ker\,S^*_z=K_b$.

Suppose that $\a>0$.
Then $h(w)^{-1}\in H^\i(\G_w)$ and $\|T^*_{h^{-1}}\|=\|h^{-1}\|_\i=\a^{-1}$.
Let $F\in N_\varphi\ominus K_b$.
We can write $(L(0)F)(w)=b(w)f(w)$. Then by (2.1)
\begin{eqnarray*}
\|F\|^2&=&\Big\|\sum^\infty_{n=0}z^nT^{*n}_\varphi bf\Big\|^2\\
&=& \sum^\infty_{n=0}\|T^{*n}_\varphi bf\|^2\\
&\ge& \|f\|^2+\|T^*_\varphi bf\|^2\\
&=& \|f\|^2+\|T^*_h f\|^2\\
&=& \|f\|^2+\a^2\a^{-2}\|T^*_h f\|^2\\
&=& \|f\|^2+\a^2\|T^*_{h^{-1}}\|^2\|T^*_h f\|^2\\
&\ge& \|f\|^2+\a^2\|f\|^2\qquad\text{by Lemma \ref{lemma5.3}}\\
&=&(1+\a^2)\|L(0)F\|^2.
\end{eqnarray*}
Since by Lemma 2.2 $\|S^*_zF\|^2+\|L(0)F\|^2=\|F\|^2$,
$$
\|S^*_zF\|^2=\|F\|^2-\|L(0)F\|^2\ge\Big(1-\frac{1}{1+\a^2}\Big)\|F\|^2=
\frac{\a^2}{1+\a^2}\|F\|^2.
$$
This implies that $S^*_z$ is bounded below on $N_\varphi\ominus K_b$,
and hence $S^*_z$ has a closed range.

Suppose that $\a=0$.
Let $\{w_k\}_k$ be a sequence in $D$ satisfying 
$|h(w_k)|<1$ and $h(w_k)\to0$ as $k\to\infty$.
Let
$$
F_k(z,w)=\frac{b(w)}{1-\overline{w_k}w}+
\sum^\infty_{n=1}z^n\frac{\overline{b(w_k)}^{(n-1)}\,\overline{h(w_k)}^n}
{1-\overline{w_k}w}.
$$
Then 
$$
\|F_k\|^2\ge\Big\|\frac{1}{1-\overline{w_k}w}\Big\|^2.
$$
Using the fact that $T_g^*(1/(1-\bar{w_k}w))=\overline{g(w_k)}(1/(1-\bar{w_k}w))$ for every $g\in H^2(\G_w)$, we have
$$
F_k(z,w)=\sum^\infty_{n=0} z^nT^{*n}_\varphi
\frac{b(w)}{1-\overline{w_k}w}\in N_\varphi\ominus K_b,
$$
and therefore $$
S^*_zF_k=\sum^\infty_{n=0}z^n\frac{\overline{b(w_k)}^n\,
\overline{h(w_k)}^{(n+1)}}
{1-\overline{w_k}w},
$$
and
$$
\|S^*_zF_k\|^2\le \Big\|\frac{1}{1-\overline{w_k}w}\Big\|^2
\frac{|h(w_k)|^2}{1-|h(w_k)|^2}.
$$
It follows
$$
\|S^*_zF_k\|^2\le \frac{|h(w_k)|^2}{1-|h(w_k)|^2}\|F_k\|^2.
$$
This implies that $S^*_z$ is not bounded below on $N_\varphi\ominus K_b$.
Since $S^*_z$ is one-to-one on $N_\varphi\ominus K_b$,
$S^*_z(N_\varphi\ominus K_b)$ is not a closed subspace.
Since $S^*_z(N_\varphi)=S^*_z(N_\varphi\ominus K_q)$,
$S^*_z$ does not have a closed range.

Next we shall prove that $S^*_zN_\varphi=N_\varphi$ when $\alpha>0$.
Let $g(w)\in L(0)N_\varphi$. We have
\begin{eqnarray*}
\sum^\infty_{n=0} \|T^{*n}_\varphi T^*_{h^{-1}}bg\|^2
&=&\|T^*_{h^{-1}}bg\|^2 
+ \sum^\infty_{n=1} \|T^{*(n-1)}_\varphi g\|^2\\
&\le& \|h^{-1}\|_\infty^2\|g\|^2+\|L(0)^{-1}g\|^2\\
&<&\infty.
\end{eqnarray*}
Hence $T^*_{h^{-1}}bg\in L(0)N_\varphi$, and
\begin{eqnarray*}
S^*_zL^{-1}_0T^*_{h^{-1}}bg&=&\sum^\infty_{n=1}z^{n-1}T^{*n}_\varphi 
T^*_{h^{-1}}bg\\
&=& \sum^\infty_{n=1}z^{n-1}T^{*(n-1)}_\varphi g\\
&=&L^{-1}_0 g.
\end{eqnarray*}
This implies that $S^*_zN_\varphi=N_\varphi$.
\end{proof}

\begin{corollary}\label{corollary5.5}
With notations as in Theorem 3.5, the following conditions are equivalent.
\begin{itemize}
\item[(i)]
$\a\not=0$.
\item[(ii)]
$S^*_z$ has a closed range.
\item[(iii)]
$S^*_zN_\varphi=N_\varphi$.
\item[(iv)]
$T^*_\p L(0)N_\p=L(0)N_\p$.
\end{itemize}
\end{corollary}

Theorem 3.5 in particular shows that $S_z$ is injective when $\alpha>0$. This is in fact a general phenomenon 
on $N_\p$. The following fact (cf. \cite[p.85]{Gar}) is need to this end.

\begin{lemma}\label{lemma5.6}
Let $h(w)$ be an outer function on $\Gamma_w$.
Then there is a sequence of outer functions
$\{h_k\}_k$ in $H^\infty(\Gamma_w)$ such that
$\|h_kh\|_\infty\le 1$ and $h_kh\to 1$ a.e. on $\Gamma_w$
as $k\to \infty$.
\end{lemma}

\begin{theorem}\label{theorem5.7}
$S_z$ is injective on $N_{\p}$.
\end{theorem}

\begin{proof}
We show that $S^*_z$ has a dense range. Let $\varphi(w)=b(w)h(w)$ be the inner-outer factorization of $\varphi$.
By Lemma \ref{lemma5.6}, there is a sequence $\{h_k\}_k$ in $H^\infty(\Gamma_w)$ such that
$$
\|h_kh\|_\infty\le 1\,\,\,\text{and}\,\,\,\, h_kh\to 1\,\,\,\,
\text{a.e. on $\Gamma_w$ as $k\to \infty$.}\leqno{(3.3)}
$$
Let $g(w)\in L(0)N_\varphi$.
By Lemma \ref{lemma5.3}, we have
\begin{eqnarray*}
\sum^\infty_{n=0}\|T^{*n}_\varphi T^*_{h_k}bg\|^2
&=&\|T^*_{h_k}bg\|^2+\sum^\infty_{n=1}\|T^*_{h_kh}T^{*(n-1)}_\varphi g\|^2\\
&\le& \|h_k\|^2_\infty\|g\|^2+\sum^\infty_{n=1}\|T^{*(n-1)}_\varphi g\|^2
\qquad\text{by (3.3)}\\
&=& \|h_k\|^2_\infty\|g\|^2+\|L(0)^{-1}g\|^2\\
&<&\infty.
\end{eqnarray*}
Hence $T^*_{h_k}bg\in L(0)N_\varphi$, and we have
\begin{eqnarray*}
&&\|S^*_zL(0)^{-1}T^*_{h_k}bg-L(0)^{-1}g\|^2\\
&=& \sum^\infty_{n=0}\|T^{*(n+1)}_\varphi T^*_{h_k}bg
-T^{*n}_\varphi g\|^2\\
&=&\sum^\infty_{n=0}\|T^*_{h_kh-1}T^{*n}_\varphi g\|^2\\
&\le& \sum^\infty_{n=0}\|(\overline{h_kh}-1)T^{*n}_\varphi g\|^2\\
&=& \int^{2\pi}_0 |(hh_k)(e^{i\theta})-1|^2\sum^\infty_{n=0}|
(T^{*n}_\varphi g)(e^{i\theta})|^2\,\frac{d\theta}{2\pi}.
\end{eqnarray*}
Since $g\in L(0)N_\varphi$, 
$$
\sum^\infty_{n=0}|T^{*n}_\varphi g|^2\in L^1(\Gamma_w).
$$
Hence by (3.3) and the Lebesgue dominated convergence theorem,
$$
\|S^*_zL^{-1}_0T^*_{h_k}bg-L^{-1}_0 g\|^2\to0\quad\text{as $k\to\infty$.}
$$
This implies that $S^*_z$ has a dense range.
\end{proof}

\begin{corollary}\label{corollary5.11}
Let $\p(w)=b(w)h(w)$ be the inner-outer factorization of $\p(w)$.
Then the following are equivalent.
\begin{itemize}
\item[(i)]
$S_z$ is Fredholm.
\item[(ii)]
$b(w)$ is a finite Blaschke product and $h^{-1}(w)\in H^\i(\G_w)$.
\end{itemize}
In this case, 
$-\,ind\,(S_z)$ is the number of zeros of $b(w)$ in $D$ counting multiplicites.
\end{corollary}
\begin{proof}
We let $\a=\displaystyle{\inf_{w\in D}}|h(w)|$.
$S_z$ is Fredholm if and only if $S^*_z$ is Fredholm, and 
by Lemma 3.4 and Theorem 3.5 this is equivalent to $b$ being a 
finite Blaschke product and $\a>0$. Clearly, $\a>0$ if and only if $h^{-1}(w)\in H^\i(\G_w)$. 

\end{proof}

A quotient module $N$ is said to be {\em essentially reductive} if both $S_z$ and $S_w$ are essentially normal, i.e.,
$[S^*_z, S_z]$ and $[S^*_w, S_w]$ are both compact. Essential reductivity is an important concept and has been studied 
recently in various contexts. In the context here, it will be interesting to see what type of $\p$ makes 
$N_\p$ essentially reductive. Proposition 3.2 has a couple of consequences to this end. 
A general study will be made in a different paper.

\begin{corollary}
For every $\p\in H^2(\G_w)$, $[S^*_z, S_w]$ is Hilbert-Schmidt on $N_\p$.
\end{corollary}
\begin{proof}
We let $R_z$ and $R_w$ denote the multiplications by $z$ and $w$ on the submodule $M_\p$, respectively.
It then follows from Proposition 3.2 and Theorem 2.3 in \cite{Y3} 
that $[R_z^*, R_z][R_w^*, R_w]$ is Hilbert-Schmidt,
and the corollary thus follows from Theorem 2.6 in \cite{Y3}.
\end{proof}

In the case $\p$ is in the disk algebra $A(D)$, there is a sequence of polynomials $p_n\rightarrow \p$ in $A(D)$, and hence
$[S^*_z, p_n(S_w)]\rightarrow [S^*_z, \p(S_w)]$ in operator norm. Since $S_z=\p(S_w)$ on $N_\p$, we easily obtain the following 
corollary.
\begin{corollary}
If $\p\in A(D)$, then $S_z$ is essentially normal.
\end{corollary}

{\noindent}{\bf Question 1.} For what $\p\in H^2(\G_w)$ is $S_w$ essentially normal on $N_\p$? \\

In the case $\p$ is inner, this question can be settled by direct calculations. We will do it in Section 5.



\section{Compactness of $L(0)|_N$ and $D_z$}
In view of Lemma 2.2, the compactness of $L(0)|_N$ or $D_z$ will give us much information about the operator 
$S_z$. So to determine whether $L(0)|_N$ or $D_z$ is compact for a certain quotient module $N$ is of great interests.
In the case of $N_\p$, the compactness is undoubtly dependent on the properties of $\p$. This section aims to unveil the
connection. 

We first look at the compactness of $L(0)|_{N\p}$. For each fixed $\zeta\in D$, we denote by $Z_\varphi(\zeta)$ 
the number of zeros of $\zeta-\varphi(w)$ in $D$ counting multiplicities. This integer-valued function 
has an important role to play in this study. As a matter of fact, 
in \cite[Theorem 5.2.2]{Y2}, the second author showed that if 
$L(0)$ on $N_\varphi$ is compact, then $Z_\varphi(\zeta)$ is a finite constant on $D$. The following study describes
the function $\p$ for which this is the case.

\begin{lemma}
Let $\varphi(w)=b(w)h(w)$ be the inner-outer factorization of $\varphi$. 
Then $Z_\varphi(\zeta)$ is a finite constant on $D$ if and only if 
$b$ is a finite Blaschke product and $|h(w)|\ge 1$ for every $w\in D$.
\end{lemma}
\begin{proof}
It is easy to see that that $b$ is a finite Blaschke product and $|h(w)|\ge 1$ for every $w\in D$ if and only if
\[\displaystyle{\liminf_{|w|\to1}}\,|\varphi(w)|\ge1.\]
Suppose that $c = Z_\varphi(\zeta)$ for every $\zeta\in D$.
To prove the necessity by contradiction, we assume that there exists 
a sequence $\{w_n\}_n$ in $D$ such that 
$\displaystyle{\sup_n}\, |\varphi(w_n)|<1$ and $|w_n|\to1$.
We may assume that $\varphi(w_n)\to\zeta_0\in D$.
Then there exists $r_0,0<r_0<1$, such that the number of zeros 
of $\zeta_0-\varphi(w)$ in $r_0D$ equals to $c$.
By the Hurwitz theorem, for a large positive integer $n_0$,
the number of zeros of $\varphi(w_{n_0})-\varphi(w)$ in $r_0D$
equals to $\alpha$. Further, we may assume that $w_{n_0}\notin r_0D$.
Hence the number of zeros of $\varphi(w_{n_0})-\varphi(w)$ in $D$ 
is greater than $c$ which contradicts the fact that $Z_\varphi(\zeta)$ is a constant.

The sufficiency is an easy consequence of Rouch\'e's theorem in Complex Analysis. In fact, 
if $b(w)$ is a finite Blaschke product and $h(w)$ is an outer function with $|h(w)|\ge1$ on $D$,
then by Rouch\'e's theorem, for each $\zeta\in D$ the number of zeros of $\zeta-\p(w)$
in $D$ coincides with the number of zeros of $b(w)$ in $D$. So $Z_\varphi(\zeta)$ is a finite constant.
\end{proof}

\begin{theorem}\label{theorem6.1}
Let $\varphi(w)=b(w)h(w)$ be the inner-outer factorization of $\varphi$.
Then the following conditions are equivalent.
\begin{itemize}
\item[(i)]
$L(0)$ on $N_\varphi$ is compact.
\item[(ii)]
$b$ is a finite Blaschke product and $|h(w)|\ge 1$ for every $w\in D$.
\end{itemize}
\end{theorem}

\begin{proof}
(i) $\Rightarrow$ (ii) If $L(0)$ on $N_\varphi$ is compact, then by Theorem 5.2.2 in \cite{Y2} $Z_\varphi(\zeta)$ is
a finite constant, and (ii) thus follows from Lemma 4.1.

(ii) $\Rightarrow$ (i) For any positive integer $m$,
we have
$$
H^2(\Gamma_w)\ominus b^m(w)H^2(\Gamma_w)=\sum^{m-1}_{j=0}
\oplus \,b^j(w)\big( H^2(\Gamma_w)\ominus b(w)H^2(\Gamma_w)\big).
$$
Since $b$ is a finite Blaschke product,
$dim\,\big( H^2(\Gamma_w)\ominus b(w)H^2(\Gamma_w)\big)<\i$ and
$H^2(\Gamma_w)\ominus b(w)H^2(\Gamma_w)$ is contained in the 
disk algebra $A(D)$. One easily sees that
$$
T^*_\p b^j(w)\big( H^2(\Gamma_w)\ominus b(w)H^2(\Gamma_w)\big)
\subset b^{j-1}(w)\big( H^2(\Gamma_w)\ominus b(w)H^2(\Gamma_w)\big),
$$
so that
$$
H^2(\Gamma_w)\ominus b^m(w)H^2(\Gamma_w)\subset L(0)N_\varphi.
$$ 
Then 
$$
L(0)N_\varphi=(H^2(\Gamma_w)\ominus b^mH^2(\Gamma_w))\oplus
(b^mH^2(\Gamma_w)\cap L(0)N_\varphi)
$$ 
and hence
$$
N_\varphi=L(0)^{-1}(H^2(\Gamma_w)\ominus b^mH^2(\Gamma_w))+
L(0)^{-1}(b^mH^2(\Gamma_w)\cap L(0)N_\p),
$$ 
which is in fact a direct sum because $L(0)|_{N_{\p}}$ is injective. 
For simplicity we write this decomposition as \[N_\varphi=N_{1,m}+N_{2,m}.\]

Since $dim\,(N_{1,m})<\infty$,
to prove that $L(0)$ on $N_\varphi$ is compact it is sufficient to prove
that $\lim_{m\to \infty}\|L(0)|_{N_{2,m}}\|=0$, i.e.,
$$
\sup_{b^mg\in L(0)N_\varphi} \frac{\|b^mg\|^2}{\|L(0)^{-1}b^mg\|^2}\to 0\quad\text{as $m\to\infty$.}
$$
Let $b^mg\in L(0)N_\varphi$ and $0\le n\le m$.
By Lemma \ref{lemma5.3}, $T^*_hb^{m-1}g=T^*_\varphi b^mg\in H^2(\Gamma_w)$, so that
$$
T^{*2}_hb^{m-2}g=T^*_hT^*_hT^*_bb^{m-1}g=
T^*_hT^*_bT^*_hb^{m-1}g=T^{*2}_\varphi b^mg\in H^2(\Gamma_w).
$$
Repeating this, we have
$$
T^{*n}_hb^{m-n}g=T^{*n}_\varphi b^mg\in H^2(\Gamma_w).\leqno{(4.1)}
$$
Using the fact that $L(0)A_\p f=f$, i.e., \[L^{-1}(0)f=\sum_{j=0}^{\infty}z^jT^{*j}_{\p}f,\] and that $\|h^{-1}\|_{\infty}\leq 1$,
we calculate that
\begin{eqnarray*}
\sup_{b^mg\in L(0)N_\varphi}\frac{\|b^mg\|^2}{\|L(0)^{-1}b^mg\|^2}
&=&\sup_{b^mg\in L(0)N_\varphi}\frac{\|g\|^2}
{\sum^\infty_{j=0}\|T^{*j}_\varphi b^mg\|^2}\\
&\le&\sup_{b^mg\in L(0)N_\varphi}\frac{\|g\|^2}
{\sum^m_{j=0}\|T^{*j}_\varphi b^mg\|^2}\\
&=&\sup_{b^mg\in L(0)N_\varphi}\frac{\|g\|^2}
{\sum^m_{j=0}\|T^{*j}_hb^{m-j}g\|^2}\qquad\text{by (4.1)}\\
&\le&\sup_{b^mg\in L(0)N_\varphi}\frac{\|g\|^2}
{\sum^m_{j=0}\|T^{*j}_{h^{-1}}\|^2\|T^{*j}_hb^{m-j}g\|^2}\\
&\le&\sup_{b^mg\in L(0)N_\varphi}\frac{\|g\|^2}
{\sum^m_{n=0}\|b^{m-n}g\|^2}\qquad\text{by Lemma \ref{lemma5.3}}\\
&=& \frac{1}{m+1}.
\end{eqnarray*}
So it follows that $\lim_{m\to \infty}\|L(0)|_{N_{2,m}}\|=0$ and  
this completes the proof.
\end{proof}

\begin{corollary}
If $L(0)$ and $R(0)$ are both compact on $N_{\p}$ then $\p$ is a finite Blaschke product.
\end{corollary}
\begin{proof}
If $R(0)$ is compact on $N_{\p}$, then by the parallel statement of Theorem 5.2.2 in \cite{Y2} for $R(0)$,
the number of zeros of $z-\p(\l)$ in $D$ is a constant with respect to $\l\in D$. Since $N_\p$ is  non-trivial, this constant
is equal to $1$. So $\|\p\|_{\infty}\leq 1$, and it follows that $\|h\|_{\infty}\leq 1$. If 
$L(0)$ is also compact on $N_{\p}$, then by Theorem 4.2 $h$ is a constant of modulous $1$, hence
$\p$ is a finite Blaschke product.
\end{proof}

In fact the converse of Corollary 4.3 is also true and we will see 
it in Section 5.\\

Next we study the compactness of $D_z$. In fact, the compactness of $D_z$ and that of $L(0)|_{N_\p}$ are
closely related. 
\begin{theorem}
If $\p$ is bounded, then $L(0)|_{N_\p}$ is compact if and only if $D_z$ is compact.
\end{theorem}
\begin{proof}
The fact that the compactness of $L(0)|_{N_\p}$ implies the compactness of $D_z$ follows from Theorem 3.7 and \cite[Theorem 5.3.1]{Y2}.

To show that the compactness of $D_z$ implies that of $L(0)|_{N_\p}$, we first check that $S_z$ is Fredholm in this case. 
If $D_z$ is compact, then by Lemma 2.2 $S_z^*S_z$ is Fredholm, 
and hence $S_z^*$ has closed range. Moreover, it follows from Theorem 3.8 that $S_z^*$ is in fact onto. So it remains to show 
that $S^*_z$ has a finite dimensional kernel. If we let $\p=bh$ be the inner-outer factorization of $\p$, then by Lemma 3.4
we need to show that $H^2(\G_w)\ominus bH^2(\G_w)$ is a finite dimensional subspace in $N_{\p}$, or equivalently, $b$ is a Blaschke 
product. For every $f\in H^2(\G_w)\ominus bH^2(\G_w)$ and integers $i, j\geq 0$, one checks that
\[\langle D_z^*f, (z-\p)z^iw^j\rangle=\langle zf, (z-\p)z^iw^j\rangle
=\langle f, z^iw^j\rangle.\]
So $D_z^*f$ is orthogonal to $(z-\p)z^iw^j$ when $i\geq 1$. Therefore,
\begin{align*}
\|D^*_zf\|&=\|P_{M_{\p}}zf\|\\
&\ge\sup_{\|(z-\p)p\|\leq 1}|\langle zf, (z-\p)p\rangle|,\ \ 
\text{p are polynomials in $H^2(\G_w)$}\\
&=\sup_{\|(z-\p)p\|\leq 1}|\langle f, p\rangle|.
\end{align*}
Since 
\[\|(z-\p)p\|^2=\|p\|^2+\|\p p\|^2
\leq \|p\|^2(1+\|\p\|_{\infty}^2),\] 
we have
\[\|D^*_zf\|\geq \sup_{\|p\|\leq (1+\|\p\|_{\infty}^2)^{-1/2}}|
\langle f, p\rangle|=(1+\|\p\|_{\infty}^2)^{-1/2}\|f\|,\]
which means $D_z^*$ is bounded below by a positive constant on $H^2(\G_w)\ominus bH^2(\G_w)$. Since $D_z$ is compact,
$H^2(\G_w)\ominus bH^2(\G_w)$ is finite dimensional, and this concludes that $S_z$ is Fredholm.

Now we show that $L(0)|_{N_\p}$ is compact. For this matter, we recall the equality (cf. Proposition 5.1.1 in \cite{Y2})
\[S_zD_z+(L(0)|_N)^*(L(0)|_{M\ominus zM})=0.\]
Since $D_z$ is compact, $(L(0)|_{N}^*(L(0)|_{M\ominus zM})$ is compact. Since we have shown that $S_z$ is Fredholm in this case, 
$L(0)|_{M_\p\ominus zM_\p}$ is Fredholm by Lemma 2.1, and therefore $L(0)|_{N_\p}$ is compact.
\end{proof}

The following example gives a simple illustration for the compactness of $L(0)|_{N_\p}$.\\ 

\noindent{\bf Example 1.}
We consider a function $\varphi(w)=aw$, where $a\in \mathbb{C}$
and $a\not=0$. Let
$$
R_j=\sqrt{1+|a|^2+\cdots+|a|^{2j}}
$$
and
$$
e_j=\frac{w^j+(\overline{a}z)w^{j-1}+\cdots+(\overline{a}z)^j}{R_j}.
$$
Then it is not difficult to check that $\{e_j\}_j$ is an orthonormal basis of $N_\varphi$, and one verifies that
\[\|L(0)e_j\|^2=\Big\|\frac{w^j}{R_j}\Big\|^2=R_j^{-2}.\]
So if $|a|<1$, then $\|L(0)e_j\|^2\geq 1-|a|^2$ and hence $L(0)$ on $N$ is not compact.
If $|a|\geq 1$, then $\lim_{j\to \infty} \|L(0)e_j\|=0$ which shows that $L(0)$ on $N$ is compact. 

It is clear by Corollary 3.11 that $S_z$ is essentially normal in this case. It is easy to give a direct calculation of 
$[S^*_z,S_z]$. In fact,
$$
S_ze_j=\frac{aR_j}{R_{j+1}}e_{j+1},\quad
S^*_ze_j=\frac{\overline{a}R_{j-1}}{R_j}e_{j-1},
$$
so
\begin{eqnarray*}
&&(S^*_zS_z-S_zS^*_z)e_j\\
&=&|a|^2\Big(
\frac{R^2_j}{R^2_{j+1}}-
\frac{R^2_{j-1}}{R^2_j}\Big)e_j\\
&=&\Big(
\frac{|a|^2+\cdots +|a|^{2(j+1)}}{1+|a|^2+\cdots +|a|^{2(j+1)}}-
\frac{|a|^2+\cdots +|a|^{2j}}{1+|a|^2+\cdots +|a|^{2j}}\Big)e_j\\
&:=&c_je_j.
\end{eqnarray*}
It is clear that $c_j\to0$ as $j\to\infty$. One also observes that $S_z$ 
on $N_{aw}$ is hyponormal.

\vspace{4mm}

By \cite{IY}, we know that $\|S_z\|=\|\p\|_\i$ if $\|\p\|_\i\le1$,
and $\|S_z\|=1$ for other cases.
In the last part of this section, we calculate the norm and the essential 
norm of $L(0)|N_\p$ and $S_z$. First we recall that the essential norm $\|A\|_e$ is the norm of $A$ in the 
Calkin algebra.

Since $\|S^*_zF\|^2+\|L(0)F\|^2=\|F\|^2$ for every $F\in N_\p$, we have
$$
\|S^*_z\|^2=\sup_{F\in N_\p,\|F\|=1}\|S^*_zF\|^2=1-\inf_{F\in N_\p,\|F\|=1}
\|L(0)F\|^2,
$$
$$
\inf_{F\in N_\p,\|F\|=1}\|S^*_zF\|^2=1-\sup_{F\in N_\p,\|F\|=1}
\|L(0)F\|^2=1-\|L(0)\|^2.\leqno{(4.2)}
$$
Hence
$$
\inf_{F\in N_\p,\|F\|=1}\|L(0)F\|=
\left\{
\begin{array}{rl}
\sqrt{1-\|\p\|^2_\i},&\quad \mbox{if $\|\p\|_\i\le1$}\\
0,&\quad\mbox{other cases.}
\end{array}\right.
$$

\begin{proposition}
Let $\alpha=\displaystyle{\inf_{w\in D}}|\varphi(w)|$.
Then $\alpha < 1$ and $\|L(0)|N_\p\|=\sqrt{1-\alpha^2}$.
\end{proposition}

\begin{proof}
By \cite[Corollary 2.7]{IY}, $\varphi(D)\cap D\not=\emptyset$.
Hence $\alpha<1$. Let
$$
F=\frac{2}{(1-\overline{\p(w_0)}z)(1-\overline{w}_0w)}
$$
Let $w_0\in \Omega_\p$.
Then by Lemma 2.5, $F\in N_\p$ and
$$
\frac{\|L(0)F\|^2}{\|F\|^2}=1-|\varphi(w_0)|^2.
$$
This implies $1-|\varphi(w_0)|^2\le\|L(0)\|^2$.
Thus we get
$$
\sqrt{1-\alpha^2}\le\|L(0)\|\le1.\leqno{(4.3)}
$$
If $\alpha=0$, then $\|L(0)\|=1$.

Suppose that $\alpha>0$.
Then $(1/\varphi)(w)\in H^\infty(\Gamma_w)$, and by 
Lemma 2.4 we have
$T^*_{1/\varphi^n}T^{*n}_\varphi=I$ on $L(0)N_\varphi$ 
for every $n\ge1$.
Let $h\in L(0)N_\varphi$. We have
\begin{eqnarray*}
\|h\| &=& \|T^*_{1/\varphi^n}T^{*n}_{\varphi}h\|\\
&\le& \|T^*_{1/\varphi^n}\|\|T^{*n}_{\varphi}h\|\\
&=& \|1/\varphi\|^n_\infty\|T^{*n}_{\varphi}h\|\\
&=& \|T^{*n}_{\varphi}h\|/\alpha^n.
\end{eqnarray*}
Then $\alpha^n\|h\|\le\|T^{*n}_{\varphi}h\|$
for every $h\in L(0)N_\varphi$ and $n$.
Hence
$$
\|h\|^2\frac{1}{1-\alpha^2}\le 
\sum^\infty_{n=0}\|T^{*n}_\varphi h\|^2=\|L^{-1}_0h\|^2
$$
for every $h\in L(0)N_\varphi$,
and $\|L(0)F\|^2\le (1-\alpha^2)\|F\|$ for every $F\in N_\varphi$.
Therefore $\|L(0)\|\le\sqrt{1-\alpha^2}$.
By (4.3), $\|L(0)\|=\sqrt{1-\alpha^2}$.
\end{proof}

A combination of (4.2), Propositions 3.1 and Proposition 4.5 leads to the following

\begin{corollary}
Let $\a=\displaystyle{\inf_{w\in D}}\,|\varphi(w)|$.
Then $S_z^*$ is invertible if and only if $\alpha >0$.
In this case, 
$$
\|S^{*-1}_z\|^{-1}=\inf_{F\in N_\p,\|F\|=1} \|S^*_zF\|=\a.
$$
\end{corollary}

\begin{theorem}
Let $\varphi(w)\in H^\i(\Gamma_w)$ with $N_\varphi\not=\{0\}$.
Let $\varphi(w)=b(w)h(w)$ be the outer-inner factorization of $\varphi$.
Suppose that $L(0)$ on $N_\varphi$ is not compact.
Let $\gamma=\displaystyle{\liminf_{|w|\to1}} |\varphi(w)|$.
Then $\gamma<1$ and $\|L(0)\|_e=\sqrt{1-\gamma^2}$.
Moreover $\|L(0)\|_e\not=\|L(0)\|$ if and only if $b(w)$ is a non-constant finite Blaschke product and $1/h(w)\in H^\infty(\Gamma_w)$.
\end{theorem}

\begin{proof}
By Theorem 4.2, $\gamma<1$.
Take a sequence $\{w_j\}_j$ in $D$ such that $|\varphi(w_j)|\to\gamma$ and $|w_j|\to1$ as $j\to\infty$.
We have
\begin{eqnarray*}
\|L(0)k_{w_j}\|&=&\sqrt{1-|w_j|^2}\sqrt{1-|\varphi(w_j)|^2}\Big\|\frac{1}
{1-\overline{w}_0w}\Big\|\\
&=&\sqrt{1-|\varphi(w_j)|^2}\\
&\to&\sqrt{1-\gamma^2}.
\end{eqnarray*}
Let $K$ be a compact operator from $N_\varphi$ to $H^2(\Gamma_w)$.
Since $k_{w_j}\to0$ weakly in $N_\p$,
$\|(L(0)+K)k_{w_j}\|\to\sqrt{1-\gamma^2}$.
Hence $\|L(0)\|_e\ge\sqrt{1-\gamma^2}$.

Suppose that $\gamma=0$.
Then $1\le\|L(0)\|_e\le\|L(0)\|\le1$.
In this case, either $b$ is not a finite Blaschke product or $1/h\notin H^\infty(\Gamma_w)$.

Suppose that $0<\gamma<1$.
Then $b$ is a finite Blaschke product.
By Proposition 4.5, $\|L(0)\|=\sqrt{1-\alpha^2}$, where 
$\alpha=\displaystyle{\inf_{w\in D}}|\varphi(w)|$.
We note that $\alpha\le\gamma$.
If $\alpha=\gamma$, then we have $\|L(0)\|=\|L(0)\|_e$.
In this case, $b$ is a constant function and $1/h\in H^\infty(\Gamma_w)$.

If $\alpha<\gamma$, then $b$ is a non-constant finite Blaschke product and $1/h \in H^\infty(\Gamma_w)$.
This implies that $\alpha=0$ and $\|L(0)\|=1$.
In this case we shall prove that $\|L(0)\|_e=\sqrt{1-\gamma^2}$.
We note that $\|1/h\|_\infty=1/\gamma$.
The idea of the proof is the same as that of Theorem 4.2. We have
\begin{eqnarray*}
\sup_{b^mg\in L(0)N_\varphi}\frac{\|b^mg\|^2}{\|L^{-1}(0)b^mg\|^2}
&\le& \sup_{b^mg\in L(0)N_\varphi}\frac{\|g\|^2}
{\sum^m_{n=0}\|T^{*n}_hb^{m-n}g\|^2}\\
&=& \sup_{b^mg\in L(0)N_\varphi}\frac{\|g\|^2}
{\sum^m_{n=0}\gamma^{2n}\|T^{*n}_{1/h}\|^2\|T^{*n}_hb^{m-n}g\|^2}\\
&\le& \frac{1}{\sum^m_{n=0}\gamma^{2n}}.
\end{eqnarray*}
Hence $\|L(0)\|_e\le\sqrt{1-\gamma^2}$, so that we obtain 
$\|L(0)\|_e=\sqrt{1-\gamma^2}$.
\end{proof}

\begin{theorem}
Let $\varphi(w)\in H^2(\Gamma_w)$ satisfying
and $\varphi\not=0$.
Then $\|S_z\|_e=\|S_z\|$.
\end{theorem}

\begin{proof}
First, suppose that $0<\|\p\|_\i<1$.
Let $K$ be a compact operator on $N_\varphi$.
Let $\{w_j\}_j$ be a sequence in $D$ such that 
$|\varphi(w_j)|\to\|\varphi\|_\infty$ as $j\to\infty$.
Then $Kk_{w_j} \to 0$ as $j\to\infty$.
One easily sees that $\|S^*_zk_{w_j}\|=|\p(w_j)|$, so that 
$\|S^*_zk_{w_j}\|\to\|\varphi\|_\infty$ as $j\to\infty$.
Hence $\|S^*_z+K\|\ge\|\varphi\|_\infty$.
By \cite[Proposition 3.5]{IY}, $\|S^*_z\|=\|\varphi\|_\infty$, so that
$\|S_z\|_e=\|S^*_z\|_e=\|\varphi\|_\infty=\|S_z\|$.

Next, suppose that $1\le \|\p\|_\i\le\i$.
By \cite[Proposition 3.5]{IY}, $\|S_z\|=1$.
Suppose that $\displaystyle{\liminf_{|w|\to1}}|\varphi(w)|\ge1$.
By Theorem 4.2, $L(0)$ is compact on $N_\varphi$.
Since $S_zS^*_z=I-L^*(0)L(0)$, $\|S_zS^*_z\|_e=1$,
so that $\|S_z\|_e=1$.

Suppose that $\displaystyle{\liminf_{|w|\to1}}|\varphi(w)|<1$.
Then there exists a sequence $\{\alpha_j\}_j\subset D$ 
satisfying the following conditions;
$|\alpha_j|\to1$ as $j\to\infty$, and 
for each $j$ there exists a sequence $\{w_{j,l}\}_l$
in $\Omega_\varphi$ such that $|w_{j,l}|\to1$ and 
$\varphi(w_{j,l})\to\alpha_j$ as $l\to\infty$.
Let $K$ be a compact operator on $N_\varphi$.
Then $\|(S^*_z+K)k_{w_{j,l}}\|\to|\alpha_j|$ as $l\to\infty$.
Since $|\a_j|\to1$, we have $\|S^*_z+K\|\ge1$.
Hence $\|S_z\|_e=\|S^*_z\|_e=1=\|S_z\|$.
\end{proof}

\section{the case when $\p$ is inner}

This section gives a detailed study for the case when $\p$ is inner. 
On the one hand, the fact that $\p$ is inner makes this case very computable, 
and, as a consequence, many of the earlier results have a clean illustration in this case.
On the other hand, the case has a close connection with the two classical spaces, namely
the quotient space $H^2(\G)\ominus \p H^2(\G)$ and the Bergman space $L^2_a(D)$. This fact
suggests that the space $N_\p$ indeed has very rich structure.

Some preparations are needed to start the discussion.
With every inner function $\theta (w) $ in the Hardy space $\hho$ over the unit circle 
$\Gamma$, there is an associated contraction
$S(\theta)$ on $\hho\ominus \theta \hho$ defined by
\[S(\theta)f=P_{\theta}wf,\ \ \ f\in \hho\ominus \theta \hho,\]
where $P_{\theta}$ is the projection from $\hho$ onto 
$\hho\ominus \theta \hho$. 
The operator $S(\theta)$ is the classical Jordan block,
and its properties have been very well studied (cf. \cite{Be, SF}). We will state some of the related facts later 
in the section. Here, we display an orthonormal basis for $N_\p$.

\begin{lemma}\label{theorem4.4}
Let $\p(w)$ be a one variable non-constant inner function.
Let $\{\l_k(w)\}_{k=0}^m$ be an orthonormal basis of 
$H^2(\G_w)\ominus \p(w)H^2(\G_w)$, and
$$
e_j=\frac{w^j+w^{j-1}z+\cdots+z^j}{\sqrt{j+1}}
$$
for each integer $j\ge0$.
Then 
$\{\l_k(w)e_j(z,\p(w)); k=0,1,2,\cdots,m,  j=1,2,\cdots\}$ is an othonormal basis for $N_{\p}$.
\end{lemma}

\begin{proof}
First of all, we have the facts that
$$
N_\p=\Big\{A_\p f: f\in H^2(\G_w),\sum^\i_{n=0}\|T^*_{\p^n}f\|^2<\i\Big\},
$$
and
$$
H^2(\G_w)=\sum^\i_{j=0}\oplus\,\p^j(w)\big(H^2(\G_w)\ominus \p(w)H^2(\G_w)\big).
$$
Write
$$
E_{k,j}=\l_k(w)e_j(z,\p(w)).
$$
Then if $(k,j)\not=(s,t)$ and $j\le t$,
\begin{eqnarray*}
&&\langle E_{k,j},E_{s,t}\rangle\\
&=&\frac{1}{\sqrt{j+1}\sqrt{t+1}}
\sum^j_{l=0}\sum^t_{i=0}\big\langle \l_k(w)\p^{j-l}(w)z^l,
\l_s(w)\p^{t-i}(w)z^i\big\rangle\\
&=&\frac{(j+1)\big\langle \l_k(w),
\p^{t-j}(w)\l_s(w)\big\rangle}{\sqrt{j+1}\sqrt{t+1}}\\
&=&0,
\end{eqnarray*}
and $\|E_{k,j}\|=1$ for every $k,j$.
Let $f(w)\in H^2(\G_w)$ and write
$$
f(w)=\sum^\i_{j=0}\oplus\Big(\sum^m_{k=1}a_{k,j}\l_k(w)\Big)\p^j(w),\quad \sum^\i_{j=0}\sum^m_{k=0}|a_{k,j}|^2<\i.
$$
Then
$$
\sum^\i_{n=0}\|T^*_{\p^n}f(w)\|^2=\sum^\i_{n=0}
\sum^\i_{j=n}\sum^m_{k=0}|a_{k,j}|^2=
\sum^\i_{j=0}(j+1)\sum^m_{k=0}|a_{k,j}|^2.
$$
Hence
$$
\sum^\i_{n=0}z^nT^*_{\p^n}f(w)\in N_\p \iff
\sum^\i_{j=0}(j+1)\sum^m_{k=0}|a_{k,j}|^2<\i.
$$
In this case, we have
\begin{eqnarray*}
&&\sum^\i_{n=0}z^nT^*_{\p^n}f(w)\\
&=&\sum^\i_{j=0}\Big(
\sum^m_{k=0}a_{k,j}\l_k(w)\Big)
(\p^j(w)+\p^{j-1}(w)z+\cdots+z^j)\\
&=&\sum^\i_{j=0}\sum^m_{k=0}\sqrt{j+1}a_{k,j}E_{k,j}.
\end{eqnarray*}
This shows that $\{E_{k,j}\}_{k,j}$ is an othonormal basis of 
$N_\p=H^2(\G^2)\ominus M_{\p}$.
\end{proof}

The operators $L(0)|_{N_\p}$, $R(0)|_{N_\p}$ and $D_z$ are easy to 
calculate in this case.
In fact, one checks that
$$
L(0)E_{k,j}=\frac{\l_k(w)\p^j(w)}{\sqrt{j+1}},
$$
and
$$
R(0)E_{k,j}=\frac{\l_k(0)(\p(0)^j+\p(0)^{j-1}z+\cdots+z^j)}{\sqrt{j+1}}.
$$
So $L(0)|_{N_\p}$ and $R(0)|_{N_\p}$ are both compact if $m<\i$, that is, $\p(w)$ is a finite Blaschke product. 
We summarize this observation and Corollary 4.3 in the following corollary.

\begin{corollary}
For $\p\in H^2(\G_w)$, $L(0)$ and $R(0)$ are both compact on $N_\p$ if and only if $\p$ is a finite Blaschke product.
\end{corollary}

The operator $D_z$ is also easy to calculate in this case. One first verifies that 
$$
X_{k,j}:=\frac{\l_k(w)}{\sqrt{j+2}}\big(ze_j(z,\p(w))-\sqrt{j+1}\p^{j+1}(w)\big),
\quad 0\le k\le m,\,\,\, 0\le j<\i,
$$
is an othonormal basis for $M_{\p}\ominus zM_{\p}$. Then
$$
D_zX_{k,j}=\frac{\l_k(w)e_j(z,\p(w))}
{\sqrt{j+2}}=\frac{1}{\sqrt{j+2}}E_{k,j}
$$
which is also compact if $\p(w)$ is a finite Blaschke product.
 
Two other observations are also worth mentioning. First one calculates that
\begin{eqnarray*}
\langle zE_{k,j},E_{s,t}\rangle&=&\frac{1}{\sqrt{j+1}\sqrt{t+1}}
\sum^j_{l=0}\sum^t_{i=0}
\big\langle z\l_k(w)\p^{j-l}(w)z^l,
\l_s(w)\p^{t-i}(w)z^i\big\rangle\\
&=&\frac{1}{\sqrt{j+1}\sqrt{t+1}}
\sum^j_{l=0}\sum^t_{i=0}
\big\langle \l_k(w),\l_s(w)\p^{t+l-i-j}(w)z^{i-l-1}\big\rangle.
\end{eqnarray*}
Hence
$$
\langle zE_{k,j},E_{s,t}\rangle=0\iff t=j+1\,\,\,\text{and}\,\,\, k=s,
$$
and
\begin{eqnarray*}
S_zE_{k,j}&=&\langle S_zE_{k,j},E_{k,j+1}\rangle E_{k,j+1}\\
&=&\frac{1}{\sqrt{j+1}\sqrt{j+2}}
\sum^j_{l=0}
\langle \l_k(w),\l_k(w)\rangle E_{k,j+1}\\
&=&\frac{\sqrt{j+1}}{\sqrt{j+2}}
E_{k,j+1}.
\end{eqnarray*}
This calculation reminds us of the Bergman shift $B$ on 
the Bergman space $L^2_a(D)$ with the orthonormal basis
$\{\sqrt{j+1}{\zeta^j}\}_j$. 
In fact, if we define the operator \[U: N_{\p}\longrightarrow \left(H^2(\Gamma)\ominus \p H^2(\Gamma)\right)\otimes L^2_a(D)\] by
\[U(E_{k,j})=\lambda_k(w)\sqrt{j+1}\zeta^j,\]
then $U$ is clearly a unitary operator, and one checks that
\begin{equation*}
US_z=(I\otimes B)U.\tag{5.1}
\end{equation*}
So from this view point $N_{\p}$ can be identified as 
$\left(H^2(\Gamma)\ominus \p H^2(\Gamma)\right)\otimes L^2_a(D)$.
As both $H^2(\Gamma)\ominus \p H^2$ and $L^2_a(D)$ are 
classical subjects, this observation indicates that the space 
$N_{\p}$ indeed has very rich structure.

The other observation is about the range $R(D_z)$. Let $F\in N_\p$. 
Then by Theorem 2.3,
$$
F\in T^*_z\big(M_{\p}\ominus zM_{\p}\big)\iff 
\sup_{G\in N_\p,\|G\|=1}\frac{|\langle S^*_z G, F\rangle|}
{\|L(0)G\|}<\i.
$$
Write
$$
F=\sum^m_{k=0}\sum^\i_{j=0}a_{k,j}E_{k,j},\quad \sum^m_{k=0}\sum^\i_{j=0}
|a_{k,j}|^2<\i,
$$
$$
G=\sum^m_{k=0}\sum^\i_{j=0}b_{k,j}E_{k,j},\quad \sum^m_{k=0}\sum^\i_{j=0}
|b_{k,j}|^2=1.
$$
Then
\begin{eqnarray*}
\frac{|\langle S^*_z G, F\rangle|}{\|L(0)G\|}
&=&
\frac{
\big|\big\langle 
\sum^m_{k=0}\sum^\i_{j=0}b_{k,j}E_{k,j}, 
\sum^m_{k=0}\sum^\i_{j=0}a_{k,j}S_zE_{k,j}
\big\rangle\big|
}
{
\big\|
\sum^m_{k=0}\sum^\i_{j=0}b_{k.j}\l_k(w)\p^j(w)
\big\|
}\\
&=&
\frac{
\big| 
\sum^m_{k=0}
\big\langle 
\sum^\i_{j=0}b_{k.j}E_{k,j}, 
\sum^\i_{j=0}a_{k,j}S_zE_{k,j}
\big\rangle
\big|
}
{
\sqrt{
\sum^m_{k=0}\sum^\i_{j=0}
\frac{|b_{k,j}|^2}
{\sqrt{j+1}}
}}\\
&=&
\frac{
\big| 
\sum^m_{k=0}\sum^\i_{j=0}
\frac{\sqrt{j+1}}
{\sqrt{j+2}}
b_{k,j+1}\overline{a}_{k,j}
\big|
}
{
\sqrt{
\sum^m_{k=0}\sum^\i_{j=0}
\frac{|b_{k,j}|^2}
{\sqrt{j+1}}
}
}
\end{eqnarray*}
and
$$
\sup_{G\in N_\p,\|G\|=1}\frac{|\langle S^*_z G, F\rangle|}{\|L(0)G\|}
=\sqrt{\sum^m_{k=0}\sum^\i_{j=0}(j+1)|a_{k,j}|^2}.
$$
Write $c_{k,j}=\sqrt{j+1}a_{k,j}$, then we have 
$F\in D_z\big(M_{\p}\ominus zM_{\p}\big)$ if and only if
$$
F=\sum^m_{k=0}\sum^\i_{j=0}
\frac{c_{k,j}E_{k,j}}{\sqrt{j+1}},\quad
\sum^m_{k=0}\sum^\i_{j=0}
|c_{k,j}|^2<\i.
$$
So 
$$
U(R(D_z))=\left(H^2(\Gamma)\ominus \p H^2(\Gamma)\right)\otimes H^2(\G).
$$

It follows directly from (5.1) that $S_z$ on $N_{\p}$ is essentially normal if and only if $\p$ is a finite Blaschke 
product. Now we take a look at the essential normality of $S_w$. Some facts about the space 
$H^2(\Gamma)\ominus \p H^2(\Gamma)$ need to be mentioned here. We recall that the Jordan block $S(\p)$ is defined by 
\[S(\p)g=P_\p wg, \ \ g\in H^2(\Gamma)\ominus \p H^2(\Gamma),\]
where $P_\p$ is the orthogonal projection from $H^2(\Gamma)$ onto $H^2(\Gamma)\ominus \p H^2(\Gamma)$. 
The two functions $P_\p 1$ and $P_\p \bar{w}\p$ play important roles here, and we let the operator $T_0$ on $H^2(\G)\ominus \p H^2(\G)$ be defined 
by $T_0g=\langle g,P_\p \bar{w}\p\rangle P_\p 1$. One verifies that
\[T_0^*T_0g=\|P_\p 1\|^2\langle g,P_\p \bar{w}\p\rangle P_\p \bar{w}\p,\ \ T_0T_0^*g
=\|P_\p \bar{w}\p\|^2\langle g,P_\p 1\rangle P_\p 1,\]
and 
\begin{align*}
I-S(\p)^*S(\p)=\|P_\p 1\|^{-2}T_0^*T_0,\ \ I-S(\p)S(\p)^*=
\|P_\p \bar{w}\p\|^{-2}T_0T_0^*.\tag{5.2}
\end{align*}

For every $g\in H^2(\Gamma)\ominus \p H^2(\Gamma)$,
we decompose $wg$ as \[wg=S(\p)g+(I-P_\p )wg.\] 
Using the facts that $(I-P_\p)wg=\langle wg, \p\rangle \p$,
$P_\p1=1-\overline{\p(0)}\p$ and $S_{\p}=S_z$, where $S_\p g=P_{N_\p}\p g$,
we have 
\begin{align*}
S_wge_j&=\sum_{m,n}\langle wg e_j,E_{m,n}\rangle E_{m,n}\\
&=\sum_{m,n}\Big\langle (S(\p)g) e_j+\langle wg, \p\rangle  \frac{\p P_\p1}{1-\overline{\p(0)}\p}e_j, E_{m,n}\Big\rangle E_{m,n}\\
&=(S(\p)g)e_j+ \langle wg, \p\rangle \sum_{m,n}\Big\langle 
\frac{\p P_\p1}{1-\overline{\p(0)}\p}e_j, E_{m,n}\Big\rangle E_{m,n}\\
&=(S(\p)g)e_j+ \langle g, \bar{w}\p\rangle (1-\overline{\p(0)}S_z)^{-1}S_z(P_\p1\cdot e_j).
\end{align*}

So 
\begin{equation*}
US_wU^*=S(\p)\otimes I+ T_0\otimes (1-\overline{\p(0)}B)^{-1}B.\tag{5.3}
\end{equation*}

For further discussion, we assume $\p$ is not a singular inner function, i.e., $\p$ has a zero in $D$. 
We first look at the case when $\p(0)=0$. In this case (5.3) reduces to the cleaner expression
\begin{equation*}
US_w=(S(\p)\otimes I+T_0\otimes B)U.\tag{5.4}
\end{equation*}
Using (5.4) and the fact $S^*(\p)T_0=0$, one easily verifies that
\[US^*_wS_wU^*=S(\p)^*S(\p)\otimes I+T_0^*T_0\otimes B^*B,\]
and \[US_wS^*_wU^*=S(\p)S(\p)^*\otimes I+T_0T_0^*\otimes BB^*.\]
Then by (5.2)
\begin{align*}
U[S_w^*, S_w]U^*&=(I-S(\p)S(\p)^*)\otimes I-(I-S(\p)^*S(\p))\otimes I\\
&+ T_0^*T_0\otimes B^*B-T_0T_0^*\otimes BB^*\tag{5.5}\\
&=T_0T_0^*\otimes (I-BB^*)-T_0^*T_0\otimes (I-B^*B).
\end{align*}
Since $T_0$ is of rank 1 and it is well-known that $I-BB^*$ and $I-BB^*$ are Hilbert-Schmidt, (5.5) implies that 
$[S_w^*, S_w]$ is Hilbert-Schmidt. The Hilbert-Schmidt norm of $[S_w^*, S_w]$ can be readily calculated in this case. First of 
all, $P_{N_\p}1=1$ and $P_{N_\p}\bar{w}\p=\bar{w}\p$.
Let $\lambda_k, k=0,1,2, ...,$ be an orthonormal basis of 
$H^2(\Gamma)\ominus \p H^2(\Gamma)$ and $\lambda_0=1$. Then by (5.5),
\begin{align*}
[S_w^*, S_w]\lambda_ke_j&=\frac{(T_0T_0^*\lambda_k)e_j}{j+1}-\frac{(T_0^*T_0 \lambda_k)e_j}{j+2}\\
&=\frac{\lambda_k(0)e_j}{j+1}-\frac{\langle\lambda_k, \bar{w}\p\rangle\bar{w}\p e_j}{j+2},
\end{align*}
and one calculates that 
\[\sum_{k}\|[S_w^*, S_w]\lambda_ke_j\|^2=\frac{1}{(j+1)^2}+\frac{1}{(j+2)^2}-\frac{2|\p'(0)|^2}{(j+1)(j+2)},\]
from which it follows that

\begin{equation*}
\|[S_w^*, S_w]\|_{H.S}^2=\frac{\pi^2}{3}-1-2|\p '(0)|^2.
\end{equation*}

In the case $\p(0)\neq 0$, we need an additional general fact.
For $\alpha\in D$, we let $x_{\alpha}(w)=\frac{\alpha -w}{1-\bar{\alpha}w}$.
So if we let operator $U_{\alpha}$ be defined by
 \[U_{\alpha}(f)(z,w):=\frac{\sqrt{1-|\alpha|^2}}{1-\bar{\alpha}w} f(z,x_{\alpha}(w)),\ \ f\in H^2(D^2),\]
then it is well-known that $U_{\alpha}$ is a unitary. We let $M'=U_{\alpha}([z-\p])=[z-\p(x_{\alpha})]$ and $N'=H^2(D^2)\ominus M'$.
The two variable Jordan block on $N'$ is denoted by $(S'_z, S'_w)$. 
Then by \cite{Y6}, 
\begin{equation*}
U_{\alpha}S_zU_{\alpha}^*=S'_z,\ \ U_{\alpha}S_wU_{\alpha}^*=x_{\alpha}(S'_w).
\end{equation*}
Since $x_{\alpha}(x_{\alpha}(w))=w$, we also have 
\begin{equation*}
U_{\alpha}x_{\alpha}(S_w)U_{\alpha}^*=S'_w.
\end{equation*}

So if $\p(0)\neq 0$, we pick any zero of $\p$, say $\alpha$.
Since $\p(x_a(0))=\p(\alpha)=0$, $[{S'_w}^*, S'_w]$ is Hilbert-Schmidt by the above calculations, and it then follows that
$[S_w^*, S_w]$ is Hilbert-Schmidt (cf. \cite[Lemma 1.3]{Y1}). So in conclusion, when $\p$ is not singular $[S_w^*, S_w]$ is 
Hilbert-Schmidt on $N_\p$.

These calculations on $S_z$ and $S_w$ prove the following theorem.
\begin{theorem}
Let $\p$ be an one variable inner function. Then $N_{\p}$ is essentially reductive if and only if 
$\p$ is a finite Blaschke product.
\end{theorem} 

On $N_{\p}$, the commutater $[S_z^*, S_w]$ can also be easily calculated. One sees that
\begin{align*}
US_z^*S_wU&=(I\otimes B^*)\left(S(\p)\otimes I+T_0\otimes (1-\overline{\p(0)}B)^{-1}B\right)\\
&=S(\p)\otimes B^*+T_0\otimes B^*(1-\overline{\p(0)}B)^{-1}B,
\end{align*}
and
\begin{align*}
US_wS_z^*U&=\left(S(\p)\otimes I+T_0\otimes (1-\overline{\p(0)}B)^{-1}B\right)(I\otimes B^*)\\
&=S(\p)\otimes B^*+T_0\otimes (1-\overline{\p(0)}B)^{-1}BB^*.
\end{align*}
So \[[S_z^*, S_w]=T_0\otimes [B^*, (1-\overline{\p(0)}B)^{-1}B].\]
It was shown in \cite{Zh} that 
\begin{equation*}
tr\,[f^*(B), g(B)]=\int_{D} f'(w)\overline{g'(w)}dA, \tag{5.6}
\end{equation*}
where $f$ and $g$ are analytic functions on $D$ that are continuous on $\bar{D}$ and the derivatives $f'$ and $g'$ are in $L^2_a(D)$.
Using (5.6), one easily verifies that $[B^*, (1-\overline{\p(0)}B)^{-1}B]$ is trace class
with $tr\,[B^*, (1-\overline{\p(0)}B)^{-1}B]=1$. Therefore, $[S_z^*, S_w]$ is trace class with
\begin{align*}
tr\,[S_z^*, S_w]&=tr\,T_0\cdot tr\,[B^*, (1-\overline{\p (0)}B)^{-1}B]\\
&=tr\,T_0\\
&=\overline{\p'(0)}.
\end{align*}

\vspace{4mm}

{\noindent}{\bf Example 2.} As we have remarked before that $S_z$ on $N_w$ is equivalent to the Bergman shift $B$ and $S_z=S_w$ in this 
case, and moreover $\p'=1$. So from the calculations above
\[tr\,[B^*, B]=1,\ \ \ \text{and}\ \ \ \|[B^*, B]\|^2_{H.S.}=\frac{\pi^2}{3}-3.\]


\bibliographystyle{amsplain}

\end{document}